\documentclass[psamsfonts,reqno]{amsart}

\usepackage{amsmath}
\usepackage{amsthm}
\usepackage{amsfonts}
\usepackage{amssymb}
\usepackage{enumerate}
\usepackage{url}
\usepackage[all]{xy}

\theoremstyle{plain}

\theoremstyle{definition}
\newtheorem*{thm*}{Theorem}
\newtheorem{thm}{Theorem}[section]

\newtheorem{prop}[thm]{Proposition}
\newtheorem{lem}[thm]{Lemma}

\newtheorem{cor}[thm]{Corollary}
\newtheorem{corollary}[thm]{Corollary}
\newtheorem{defn}[thm]{Definition}
\newtheorem{example}[thm]{Example}

\newtheorem{remark}[thm]{Remark}

\newcommand\supress[1]{}

\newcommand\tuple[1]{\langle #1\rangle}

\newcommand\N{{\mathbb{N}}}
\newcommand\R{{\mathbb{R}}}

\newcommand\Z{{\mathbb{Z}}}

\newcommand\Ccal{{\mathcal{C}}}

\newcommand\Fcal{{\mathcal{F}}}

\newcommand\Gammabi{\Gamma_{\sf bi}}

\newcommand\id{{\sf id}}

\newcommand\Set{{\sf Set}}

\newcommand\Coalg{{\sf Coalg}}
\newcommand\append{{:}}
\renewcommand\phi{\varphi }
\newcommand\xto[1]{\xrightarrow{#1}}

\title[Coinductive properties of Lipschitz functions on streams]{Coinductive properties of \\Lipschitz functions on streams}
\author{Jiho Kim}
\address{Department of Mathematics, 
    Indiana University,
    Rawles Hall,
    831 East 3rd St,
    Bloomington, IN 47405}
\email{jihokim@indiana.edu}
\date{\today}

\begin{document}

\begin{abstract}
A simple hierarchical structure is imposed on the set of Lipschitz
functions on streams (i.e.~sequences over a fixed alphabet set) under
the standard metric.  We prove that sets of non-expanding and
contractive functions are closed under a certain coiterative
construction.  The closure property is used to construct new final
stream coalgebras over finite alphabets.  For an example, we show that
the $2$-adic extension of the Collatz function and certain variants
yield final bitstream coalgebras.
\end{abstract}

\maketitle

\section{Introduction}

In the realm of theoretical computer science, the categorical notion
of coalgebras gives a mathematical foundation for computational
dynamics.  In the appropriate categories, the finality of coalgebras
can be construed as denotational semantics of various models of
computation such as automata \cite{Ru1}, programming languages,
recursive programming schemes \cite{MM}, and other calculi.  It also
has connection to a diverse collection of other mathematical
pursuits---the theory of non-wellfounded sets \cite{Aczel}, modal
logic \cite{kurz,Moss,pattinson,roessinger}, fractals and
self-similarity \cite{Leinster1}---to name a few.  Such connections
raise interesting questions about the extent to which the theory of
coalgebras may be useful in more ``classical'' mathematics.

According to a long tradition of children learning the 1, 2, 3's, we
learn at the very beginning---at least implicitly---that the sequence
of natural numbers and many operations on them are defined via the
principle of \emph{induction} and \emph{iteration}.  Categorically
speaking, these principles are shadows of the universal property of a
certain initial algebra.  Then we use the natural numbers to build
other sequences (i.e.~streams) of all sorts.  Although the notion of
streams is a basic one, it is ubiquitous in both mathematics and
computer science and therefore worthy of extensive study.
Analogously, the dual notions of \emph{coinduction} and
\emph{coiteration} expressed as the universal property of certain
final coalgebras lead to novel ways of expressing and understanding
definitions of streams and operations on them.

This present paper focuses particularly on stream coalgebras and
morphisms among them.  Given the standard metric on the set of
streams, we derive some coinductive closure properties on the set of
non-expanding maps and the set of distance-preserving maps. Section
\ref{sec:streams} contains the definitions and constructions necessary
for the paper.  Of particular interest is the stratification of the
set of Lipschitz functions on streams (each level of which we call
$k$-causal functions).  The stratification is achieved in several
seemingly dissimilar ways which are proven to be equivalent in Theorem
\ref{thm:k-causal-equivalence}.  We go on to show that the family of isometric embeddings on streams can also be characterized by a similar set of criteria.  Next, we define the notion of ``woven
functions'' and explore some properties which will be essential in the main result.

Section \ref{sec:coalgebra} introduces the
category theoretic notion of coalgebras and gives some clarifying examples
beyond streams.  Section \ref{sec:results} presents the coinductive
closure property on sets of stream functions which can be stated
roughly:

\begin{thm*}
A stream function coiteratively defined by a coalgebra woven from
non-expanding (resp.~distance-preserving) maps is non-expanding
(resp.~distance-preserving).  The converse holds in the case of
distance-preserving functions.
\end{thm*}

This theorem then is applied to the $2$-adic extension of the Collatz
function and its variants which figure largely in the $3x+1$ Problem.
It gives an essentially new perspective on the $3x+1$ conjugacy map as
a coalgebra isomorphism between final stream coalgebras.  The utility
of the general approach explored in this paper may be limited in terms
of resolving the $3x+1$ Problem in particular.  Nevertheless, it
identifies additional structure within a large class of Collatz-like
dynamical systems, which may shed light on these problems as a whole.

\section{Streams}
\label{sec:streams}

Let $A$ be some alphabet set (possibly infinite) and let $A^\omega$ be
the set of sequences whose components come from $A$.  Formally, these
sequences are functions from the natural numbers
$\omega=\{0,1,2,\ldots\}$ (which act as the indices) to the alphabet
$A$.  In this paper, we will refer to these sequences as
\emph{$A$-streams} and consider Lipschitz continuous functions on
them.

\subsection{$k$-causal functions}

\begin{defn}[Metric on $A^\omega$]
Given any $\sigma,\tau\in A^\omega$, define the \emph{distance between
$\sigma$ and $\tau$} to be
$$
d(\sigma,\tau)=
\begin{cases}
0 & \text{if $\sigma=\tau$}\\
2^{-i} & \text{if $\sigma\neq \tau$}
\end{cases}
$$
where $i$ is the least index such that $\sigma(i)\neq \tau(i)$.  
\end{defn}
The
metric $d$ also satisfies the ultrametric inequality,
\begin{equation}
d(\sigma,\tau) \leq \max\{ d(\sigma,\rho),  d(\rho,\tau)\},
\label{eqn:ultrametric}
\end{equation}
for any $\rho,\sigma,\tau\in A^\omega$.  If $A$ has the discrete
topology, the topology induced by this metric is the product topology
on $A^\omega$.

\begin{defn}[$k$-causal function]
\label{defn:k-causal}
Let $k\in\Z$. A function $f\colon A^\omega\to B^\omega$ is
\emph{$k$-causal} if
$$
d(f(\sigma), f(\tau)) \leq 2^{-k} d(\sigma,\tau).
$$ In other words, $f$ is $k$-causal if and only if it is Lipschitz
continuous with constant $2^{-k}$.  Furthermore, let $\Gamma_k=\{
f\colon A^\omega\to B^\omega \mid \text{ $f$ is $k$-causal}\}$.  Also
let $\Gammabi\subset\Gamma_0$ be the subset which consists of distance-preserving functions.  We will call maps in $\Gammabi$
\emph{bicausal}.
\end{defn}

\begin{example} Consider the following examples:
\begin{enumerate}[(i)]
\item For all $\ell\leq k$, $k$-causal functions are $\ell$-causal.
That is to say, $\Gamma_k\subseteq\Gamma_\ell$.
\item The map $\beta\colon \R^\omega\to \Z^\omega$ given by
$\beta(\sigma)(n) = \lfloor\sigma(n)\rfloor$ is $0$-causal.  The
function $\beta$ is neither bicausal nor $1$-causal.  (Note that while
$\beta$ is continuous, the continuity depends on $\R$ having the
discrete topology.)
\item The identity map $\id\colon A^\omega\to A^\omega$ is bicausal.
\item The tail map $t\colon A^\omega\to A^\omega$, given by
\begin{equation}
t(\sigma)(n)= \sigma(n+1),
\label{eqn:tail-def}
\end{equation}
is $(-1)$-causal but not $0$-causal.  This function is often called
the \emph{shift map} in the context of dynamical systems.  
\item For $k\geq 0$ and a word $w\in A^k$, the map $c_w\colon
A^\omega\to A^\omega$ given by
$$
c_w(\sigma)(n)=
\begin{cases}
w(n) & \text{if $n<k$}\\
\sigma(n-k) & \text{otherwise}
\end{cases}
$$ is $k$-causal. The function $c_w$ prepends a stream
$\sigma$ with a finite word $w$.  When it is convenient, we will
denote $c_w(\sigma)$ in the infix-colon notation $w\append\sigma$.
\item If $f$ is $k$-causal and $g$ is $\ell$-causal, then $g\circ f$
is $(k+\ell)$-causal.  Furthermore, if $f$ and $g$ are bicausal, then
$g\circ f$ is bicausal.
\item $0$-causal functions are non-expanding; bicausal functions are isometric embeddings.  For $k>0$, $k$-causal
functions are contraction mappings.  
\end{enumerate}
\end{example}

Theorem \ref{thm:k-causal-equivalence} enumerates several equivalent
characterizations of $k$-causal functions.  To set the stage, for
$n\geq 0$, let $\pi_n^A\colon A^{n+1}\to A^n$ be given by
$\pi_n(\sigma)(i)=\sigma(i)$ for $0\leq i< n$.  Also, let $p_n^A\colon
A^\omega \to A^n$ for $n\geq 1$ be the \emph{prefix} map, given by
$p_n(\sigma)(i)=\sigma(i)$ for $0\leq i < n$.  (Whenever possible, we
will suppress the superscripted parameter which can be deduced from
context.)  We distinguish the map $p_1\colon A^\omega\to A$ and call
it the \emph{head} function $h$ as it identifies the ``head'' of the
stream.  For each $n\geq 0$, and $\sigma,\tau\in A^\omega$, we say
$\sigma$ and $\tau$ are \emph{$n$-prefix equivalent} if and only if
$p_n(\sigma)=p_n(\tau)$.  We denote this equivalence by
$\sigma\equiv_n \tau$.  In the extreme case where $n=0$, we have
$\sigma\equiv_0\tau$ for all $\sigma,\tau\in A^\omega$.  

\begin{lem}
We have
the following simple observations.

\begin{enumerate}[(i)]
\item \label{lem:observation-1} For
any $\sigma,\tau\in A^\omega$ and $\ell\geq 0$,
$$
p_\ell(\sigma) = p_\ell(\tau)
\quad\Longleftrightarrow\quad
\sigma\equiv_\ell\tau
\quad\Longleftrightarrow\quad
d(\sigma,\tau)\leq 2^{-\ell}
$$
\item \label{lem:observation-2}For $\ell\geq 0$, $p_\ell = \pi_\ell\circ p_{\ell+1}$.
\item \label{lem:observation-3}For $\ell\geq 0$, $\sigma = p_\ell(\sigma)\append t^{(\ell)}(\sigma)$.  In particular, $\sigma = h(\sigma)\append t(\sigma)$.  (Recall the convention that $w\append\sigma = c_w(\sigma)$ for any $w\in A^\ell$ and $\sigma\in A^\omega$.  
\item \label{lem:observation-4}For any $\sigma\in A^\omega$ and $n\geq 0$, we have $\sigma(n) = h(t^{(n)}(\sigma))$.
\end{enumerate}
\label{lem:observation}
\end{lem}

\begin{thm}[$k$-causal functions]
Let $k\in\Z$ and $f\colon A^\omega\to B^\omega$.  The following are
equivalent.
\begin{enumerate}[(i)]
\item $f$ is $k$-causal.
\item For all $i,j\geq 0$ such that $k=i-j$, 
\begin{equation}
\sigma\equiv_j\tau \quad\Longrightarrow\quad f(\sigma)\equiv_i f(\tau)
\label{eqn:k-causal} 
\end{equation}
\item For $m,n\geq 0$ with $\min\{m,n\}=0$ and $k=m-n$, $f$ is the (unique)
inverse limit of a chain of maps $\{f_\ell\}_{\ell\geq 0}$ as follows:
{
$$
\xymatrix{
A^{n} \ar[d]_{f_0} & A^{n+1} \ar[d]_{f_1} \ar[l]_-{\pi_{n}} & A^{n+2} \ar[d]_{f_2} \ar[l]_-{\pi_{n+1}} & \ar[l]_-{\pi_{n+2}}\cdots & \ar[l]  A^{n+\ell} \ar[d]_{f_{\ell}} & A^{n+(\ell+1)} \ar[d]_{f_{\ell+1}} \ar[l]_-{\pi_{n+\ell}} &  \ar[l] \cdots\\
B^{m} & B^{m+1} \ar[l]^-{\pi_m} & B^{m+2} \ar[l]^-{\pi_{m+1}} & \ar[l]^-{\pi_{m+2}}\cdots &  B^{m+\ell} \ar[l] & B^{m+(\ell+1)} \ar[l]^-{\pi_{m+\ell}} &  \ar[l]\cdots
}
$$
where $f_{\ell}\circ p_{n+\ell} = p_{m+\ell} \circ f$. 
}

\end{enumerate}
\label{thm:k-causal-equivalence}
\end{thm}

\begin{proof}
For (i)$\Rightarrow$(ii), suppose $f$ is $k$-causal, and $\sigma\equiv_j\tau$ with $k+j\geq 0$.   Then, $d(\sigma,\tau)\leq 2^{-j}$.  Therefore,
$$
d(f(\sigma),f(\tau))\leq 2^{-k} d(\sigma,\tau)\leq 2^{-(k+j)}.
$$
Letting $i=k+j$, we have $f(\sigma)\equiv_i f(\tau)$.

To show (ii)$\Rightarrow$(iii), for each $\ell\geq 0$, let the
function $f_{\ell}\colon A^{n+\ell}\to B^{m+\ell}$ be given by
$$f_\ell(w)=p_{m+\ell}\left( f(w\append\sigma)\right)$$ for $w\in
A^{n+\ell}$.  First we note that $f_\ell$ is well-defined, i.e.~it
does not depend on the choice of $\sigma$.  Let $w\in A^{n+\ell}$,
then $w\append\sigma\equiv_{n+\ell} w\append\tau$ for any
$\sigma,\tau\in A^\omega$.  Since $k=(m+\ell)-(n+\ell)$, we get $f(w\append\sigma)\equiv_{m+\ell} f(w\append\tau)$ by (\ref{eqn:k-causal}).  Then by the definition of $\equiv_{m+\ell}$, we have
$p_{m+\ell}\left( f(w\append\sigma)\right)=p_{m+\ell}\left(
f(w\append\tau)\right)$, as required.

Next, we show that $f_\ell\circ \pi_{n+\ell}= \pi_{m+\ell}\circ f_{\ell+1}$.  For $v\in A^{n+\ell+1}$, let $v=wa$ where $w\in A^{n+\ell}$ and $a\in A$.  Then,
\begin{align*}
\pi_{m+\ell}(f_{\ell+1}(v))
& =  \pi_{m+\ell}(f_{\ell+1}(wa)) &&[v=wa]\\
& =  \pi_{m+\ell}(p_{m+(\ell+1)}(f(wa\append\sigma))) && [\text{def. of $f_{\ell+1}$}]\\
& =  p_{m+\ell}(f(wa\append\sigma)) &&[p_\ell = \pi_\ell\circ p_{\ell+1}]\\
& = f_\ell(w) && [\text{def. of $f_\ell$}]\\
& = f_\ell(\pi_{n+\ell}(v)) && [w=\pi_{n+\ell}(v)]
\end{align*}
We also need to show that $f_\ell\circ p_{n+\ell} = p_{m+\ell}\circ f$.  
\begin{align*}
f_\ell(p_{n+\ell}(\sigma)) 
&= p_{m+\ell}(f(p_{n+\ell}(\sigma)\append t^{(n+\ell)}(\sigma))) &&[\text{def. of $f_\ell$}]\\
&= p_{m+\ell}(f(\sigma)) && [\sigma = p_\ell(\sigma)\append t^{(\ell)}(\sigma)]
\end{align*}

Lastly, we need to verify the universal property of projective limits, namely that if there is a map $e\colon Y_A\to Y_B$ and an associated sequence of maps, $q^A_i\colon Y_A\to A^i$ and  $q^B_i\colon Y_B\to B^i$ so that the equations
\begin{equation}
q_{i} = \pi_i\circ q_{i+1} \qquad\qquad
f_i\circ q_{n+i} = q_{m+i}\circ e 
\label{eqn:univ-prop}
\end{equation}
hold for any $i$, then there exists a unique pair of maps $r^A\colon Y_A\to A^\omega$ and $r^B\colon Y_B\to B^\omega$ so that the equations  
\begin{equation}
p_j \circ r = q_j
\qquad\qquad
f\circ r = r\circ e 
\label{eqn:univ-prop-2}
\end{equation}
hold for all $j$.  To show existence, given such $e$ and $q_i$ that satisfy (\ref{eqn:univ-prop}), let $r^A$ and $r^B$ be given by 
\begin{equation}
r^*(x)(i) = q_{j}^*(x)(i)
\end{equation}
where $j>\max\{i,m,n\}$.  The index $j$ must be greater than both $m$ and $n$ so that $q_j$ is meaningful.  Also, $j$ must be greater than $i$ so that $q_j(x)$ is a word that has a least $i$ letters.  Beyond these lower bounds, it does not matter how large $j$ is, since the assumptions (\ref{eqn:univ-prop}) imply that $q_j(x)(i)=q_{j+1}(x)(i)$.  

To verify the first equation (\ref{eqn:univ-prop-2}), fix some $j$.  Then for any $i<j$ and $x\in Y$:
$$
p_j(r(x))(i) = r(x)(i) = q_j(x)(i).
$$
Therefore $p_j\circ r = q_j$. At this point, we can show uniqueness.  Suppose $r$ and $\hat r$ both satisfy $p_j \circ r = q_j$ for any $j$. If $r\neq \hat r$, then there is some $x\in Y$ and $i\geq 0$ so that $r(x)(i)\neq \hat r(x)(i)$.  But then $q_{i+1}(x)=p_{i+1}(r(x))\neq p_{i+1}(\hat r(x)) = q_{i+1}(x)$.  This is a contradiction, therefore $r$ is unique.

To verify the the second equation (\ref{eqn:univ-prop-2}), we have for $0\leq i< j$,
\begin{align*}
r(e(x))(i) 
&= q_j(e(x))(i) &&[\text{def. of $r$}]\\
&= f_{j-m}(q_{j-k}(x))(i) &&[\text{(\ref{eqn:univ-prop}), $m=k+n$}]\\
&= f_{j-m}(p_{j-k}(r(x)))(i) &&[p_\ell\circ r = q_\ell]\\
&= p_j(f(r(x)))(i) &&[f_\ell\circ p_{n+\ell} = p_{m+\ell}\circ f]\\
&= f(r(x))(i) && [i<j]
\end{align*}
This completes the portion of the proof for (ii)$\Rightarrow$(iii).

For (iii)$\Rightarrow$(i), we show the contrapositive.  Suppose $f$
is not $k$-causal.  Then there are $\sigma,\tau\in A^\omega$ so that
$$
2^{-k}d(\sigma,\tau)< d(f(\sigma),f(\tau))\leq 1.
$$
First of all, we can surmise that $\sigma\neq\tau$ due to the strict
inequality; therefore $0<d(\sigma,\tau)=2^{-i}$ for some $i\geq 0$.  Also, since $2^{-k}d(\sigma,\tau)< 1$, we can conclude that $i+k$ is positive.  We would also like to verify that $i\geq n$.  In the case where $n=0$, there is nothing to show; and in that case $n$ is positive, we have $m=0$ and therefore
$i\geq -k = n-m = n$.  In other words, $i-n$ cannot be negative.
On one hand, because $d(f(\sigma),f(\tau))> 2^{-(k+i)}$, we
get ${f(\sigma)\not\equiv_{k+i}f(\tau)}$  (i.e.~$p_{k+i}(f(\sigma))\neq p_{k+i}(f(\tau))$).  On the other hand,
$\sigma\equiv_i\tau$ (i.e.~${p_i(\sigma)=p_i(\tau)}$) since
$d(\sigma,\tau)\leq 2^{-i}$.  Applying $f_{i-n}$ to both sides of the
equation, we get $f_{i-n}(p_i(\sigma))=f_{i-n}(p_i(\tau))$.  This
demonstrates the fact that $f_{\ell}\circ p_{n+\ell}\neq p_{m+\ell}\circ f$ (where $\ell=i-n\geq 0$) as
required for $f$ to be the inverse limit.
\end{proof}

\begin{remark}
These definitions are an extension of those given in \cite{Kim1}.
What Definition \ref{defn:k-causal} named $0$-causal and $(-1)$-causal
are called \emph{causal} and \emph{subcausal}, respectively.  The
related notion of supercausal functions, however, is not exactly the
same as $1$-causal.  Unlike supercausal functions as given in
\cite{Kim1}, if $f\colon A^\omega\to B^\omega$ is $1$-causal, then
$f(\sigma)\equiv_1 f(\tau)$ for any $\sigma,\tau\in A^\omega$.  More
generally, when $k\geq 1$, the image of $k$-causal functions must have
a common $k$-prefix.  Consequently, we have a well-defined map
$D_k\colon \Gamma_k\to B^k$ given by $D_k(f) = p_k(f(\sigma)).$
\end{remark}

\begin{prop}
Let $f\colon A^\omega\to B^\omega$ be $k$-causal for $k>0$.  Then
$f=c_w\circ \hat f$ for some $w\in B^k$ and $0$-causal function $\hat
f\colon A^\omega\to B^\omega$.
\end{prop}

\begin{proof}
Let $w=D_k(f)$ and $\hat f = t^{(k)}\circ f$.  
\end{proof}

\subsection{Woven functions}

The following definition presents a way to construct nontrivial
$k$-causal maps from a set of $(k+1)$-causal maps by, in a manner of
speaking, ``weaving them together.''

\begin{defn}[Woven function]
\label{def:woven-function}
Let $\Fcal=\{f_a\}_{a\in A}$ be an $A$-indexed set of maps from
$A^\omega$ to $B^\omega$.  With a slight
abuse of notation, we can think of $\Fcal$ as a map $\Fcal\colon
A\times A^\omega\to B^\omega$ via $\Fcal(a,\sigma) = f_a(\sigma).$  We define $T_\Fcal\colon A^\omega\to
B^\omega$ by 
\begin{equation}
T_\Fcal(\sigma)=\Fcal(h(\sigma),t(\sigma))=f_{h(\sigma)}(t(\sigma)).
\label{eqn:woven-function}
\end{equation}
$T_\Fcal$ is said to be \emph{woven from $\Fcal$}.  
\end{defn}
Intuitively, for any input stream $\sigma$, the function $T_\Fcal$
gives the image of $t(\sigma)$ under a function which $h(\sigma)$
picks out from $\Fcal$.

\begin{example}[Woven function] Consider the following examples.
\begin{enumerate}[(i)]
\item If $f_a=\id$ for all $a\in A$, the function woven from
$\{f_a\}_{a\in A}$ is the tail function $t$.  
\item If $f_a=c_a$ for all $a\in A$, the function woven from
$\{f_a\}_{a\in A}$ is the identity.
\item If $A$ is finite, we can rewrite
(\ref{eqn:woven-function}) in Definition \ref{def:woven-function} as a
definition by cases.  For instance, suppose $A=\{0,1,2\}$.  For
$\Fcal=\{f_0, f_1, f_2\}$, we have
$$
T_\Fcal(\sigma) = 
\begin{cases}
f_0(t(\sigma)) & \text{if $h(\sigma)=0$}\\
f_1(t(\sigma)) & \text{if $h(\sigma)=1$}\\
f_2(t(\sigma)) & \text{if $h(\sigma)=2$}
\end{cases}
$$
\end{enumerate}
\end{example}

\begin{lem} 
A function $T\colon A^\omega\to B^\omega$ is woven from a family of
$(k+1)$-causal functions if it is $k$-causal.  The converse holds if
$k\leq 0$.
\label{lem:k-causal-is-woven}
\end{lem}

\begin{proof}
Let $T$ be $k$-causal. For each $a\in A$,
let $f_a\colon A^\omega\to B^\omega$ be a function given by
$f_a(\sigma)=T(a\append\sigma)=T(c_a(\sigma))$.  Since $c_a$ is
$1$-causal, each $f_a$ is $(k+1)$-causal.  Let $S$ be a function woven
from $\{f_a\}_{a\in A}$.  Then,
$$
S(\sigma) = f_{h(\sigma)}(t(\sigma)) = T(h(\sigma)\append t(\sigma)) = T(\sigma).
$$
Therefore, $T=S$ is woven from a family of $(k+1)$-causal functions.

Conversely, let $\Fcal=\{f_a\}_{a\in A}$ be an $A$-indexed set of
$(k+1)$-causal functions for some $k\leq 0$.  Let $\sigma,\tau\in
A^\omega$.  Then, $\sigma=a\append\sigma'$ and $\tau=b\append\tau'$
for some $a,b\in A$ and $\sigma',\tau'\in A^\omega$.  On one hand,
suppose $a\neq b$.  Then, $d(\sigma,\tau)=1$, and
$$
d(T_\Fcal(a\append\tau'),T_\Fcal(b\append\tau'))\leq 1 =
d(\sigma,\tau)\leq 2^{-k} d(\sigma,\tau).
$$ At this last inequality, we require the fact that $k\leq 0$
so that $2^{-k}\geq 1$.  On the other hand, suppose $a=b$.  Then,
\begin{align*}
d(T_\Fcal(\sigma),T_\Fcal(\tau))
&= d(T_\Fcal(a\append\sigma'),T_\Fcal(b\append\tau')) \\
&= d(T_\Fcal(a\append\sigma'),T_\Fcal(a\append\tau')) &&[a=b]\\
&=d(f_a(\sigma'), f_a(\tau')) \\
&\leq 2^{-(k+1)}d(\sigma',\tau')&& \text{[$f_a$ is $(k+1)$-causal]}\\
&=2^{-k}d(a\append\sigma',b\append\tau') && \text{[$a=b$]}\\
&= 2^{-k} d(\sigma,\tau)
\end{align*}
The calculation shows that $T_\Fcal$ is $k$-causal.
\end{proof}

\subsection{$0$-causal and bicausal functions}

The case where $k=0$ is particularly interesting because the set of
$0$-causal functions $\Gamma_0$ is closed under composition.  Furthermore, $\Gamma_0$ contains a subfamily of functions $\Gammabi$ of bicausal functions which is also closed under composition.
From Theorem \ref{thm:k-causal-equivalence} we immediately derive the following characterization of $0$-causal functions.  

\begin{cor}[$0$-causal functions, \cite{Kim1}]
Let $f\colon A^\omega\to B^\omega$ be a stream function.
The following are equivalent.
\begin{enumerate}[(i)]
\item $f$ is $0$-causal (i.e.~non-expanding).
\item For all $n\geq 0$, 
\begin{equation}
\sigma\equiv_n\tau \quad\Longrightarrow\quad f(\sigma)\equiv_n f(\tau) 
\label{eqn:0-causal}
\end{equation}
\item $f$ is the (unique) inverse limit of a chain of maps as follows:
\begin{equation}
\xymatrix{
A^{0} \ar[d]_{f_0} & A^{1} \ar[d]_{f_1} \ar[l]_-{\pi_{0}} & A^{2} \ar[d]_{f_2} \ar[l]_-{\pi_{1}} & \ar[l]_-{\pi_{2}}\cdots & \ar[l]  A^{j} \ar[d]_{f_{j}} & A^{j+1} \ar[d]_{f_{j+1}} \ar[l]_-{\pi_j} &  \ar[l] \cdots\\
B^{0} & B^{1} \ar[l]^-{\pi_0} & B^{2} \ar[l]^-{\pi_{1}} & \ar[l]^-{\pi_{2}}\cdots &  B^{j} \ar[l] & B^{j+1} \ar[l]^-{\pi_{j}} &  \ar[l]\cdots
}
\label{eqn:0-causal-solenoid}
\end{equation}
where $f_j\circ p_j = p_j \circ f$.
\end{enumerate}
\label{cor:0-causal}
\end{cor}

A similar result holds for bicausal functions with minimal changes.  Moreover, if the domain and codomain are streams over the same finite alphabet, we can strengthen the result slightly.

\begin{cor}[Bicausal functions with arbitrary alphabet]
Let $f\colon A^\omega\to B^\omega$ be a stream function.   The following are equivalent:
\begin{enumerate}[(i)]
\item $f$ is bicausal (i.e.~distance-preserving).
\item For all $n\geq 0$, 
\begin{equation}
\sigma\equiv_n\tau \quad\Longleftrightarrow\quad f(\sigma)\equiv_n f(\tau) 
\label{eqn:bicausal}
\end{equation}
\item $f$ is the (unique) inverse limit of a chain of \emph{injective} maps $\{f_j\}_{j\geq 0}$ where $f_j\circ p_j = p_j\circ f$, as arranged in (\ref{eqn:0-causal-solenoid}).
\end{enumerate}
\label{cor:bicausal}
\end{cor}

\begin{proof} 
In light of Corollary \ref{cor:0-causal}, we only need to extend the proof for the extra conclusions. 
 
For (i)$\Rightarrow$(ii), assume $f(\sigma)\equiv_n f(\tau)$.   Then, $d(f(\sigma), f(\tau)) \leq 2^{-n}$, but since $f$ is distance-preserving $d(\sigma,\tau)=d(f(\sigma),f(\tau))$.  Therefore $\sigma\equiv_n\tau$.  

For (ii)$\Rightarrow$(iii), fix some $j\geq 0$ and suppose $f_j(w)=f_j(v)$ for some $w,v\in A^j$.  Then, $p_j(f(w\append\sigma))=p_j(f(v\append\sigma'))$, or equivalently, $f(w\append\sigma)\equiv_j f(v\append\sigma')$.  Finally, we have $w\append\sigma\equiv_j v\append\sigma'$, i.e.~$w=j$, by (\ref{eqn:bicausal}).  This shows that $f_j$ is injective for any $j$.

For (iii)$\Rightarrow$(i), we already know that $f$ must be $0$-causal, therefore 
$$
d(f(\sigma),f(\tau))\leq d(\sigma,\tau).
$$
Suppose that $d(\sigma,\tau)=2^{-j}$ for some $j\in\N$ and $\sigma,\tau\in A^\omega$.   Because $d(\sigma,\tau)=2^{-j}$, we have $p_{j+1}(\sigma)\neq p_{j+1}(\tau)$ by the definition of the metric $d$.  Since $f_{j+1}$ is injective,  $f_{j+1}(p_{j+1}(\sigma))\neq f_{j+1}(p_{j+1}(\tau))$, and because of the universal property of $f$ (in particular, $f_\ell\circ p_\ell = p_\ell\circ f$), we have $p_{j+1}(f(\sigma))\neq p_{j+1}(f(\tau))$.  That is to say, $f(\sigma)$ and $f(\tau)$ first differ at an index $i<j+1$, so 
$$
d(f(\sigma),f(\tau))=2^{-i}\geq 2^{-j}=d(\sigma,\tau).
$$  
This shows that the distance must be preserved by $f$, i.e.~$f$ is bicausal. 
\end{proof}

\begin{cor}[Bicausal functions with finite alphabet, \cite{Kim1}]
Let $f\colon A^\omega\to A^\omega$ be a stream function where $A$ is finite.   The following are equivalent:
\begin{enumerate}[(i)]
\item $f$ is bicausal.
\item For all $n\geq 0$, 
$\sigma\equiv_n\tau$ if and only if $f(\sigma)\equiv_n f(\tau)$.
\item $f$ is the (unique) inverse limit of a chain of \emph{bijective} (or equivalently, \emph{surjective}) maps $\{f_j\}_{j\geq 0}$ where $f_j\circ p_j = p_j\circ f$, as arranged in (\ref{eqn:0-causal-solenoid}).
\item \label{item:causal-bijection} $f$ is a $0$-causal bijection.
\end{enumerate}
\label{cor:bicausal-finite-alphabet}
\end{cor}

\begin{proof}
If $A$ is finite, $A^j$ is finite for any finite $j$.  Consequently, $f_j\colon A^j\to A^j$ is injective if and only if it is bijective.  In light of Corollary \ref{cor:bicausal}, this shows the equivalence of (i), (ii), and (iii).

For (iv)$\Rightarrow$(iii), assume $f$ is bijective.  It is surjective, in particular, so for any word $w\in A^j$ and $\sigma\in A^\omega$, there is some $\tau\in A^\omega$ so that $f(\tau)=w\append \sigma$.  Then 
$$
f_j(p_j(\tau)) = p_j(f(\tau)) = p_j(w\append\sigma) = w
$$
That is to say, $f_j$ is surjective.  (In fact, this shows that surjective $0$-causal functions are inverse limits of surjections in general.) Since $A$ is finite, $f_j$ is also bijective.

For (iii)$\Rightarrow$(iv), we do not require that $A$ be finite.  Notice that $\{f_j^{-1}\}_{j\geq 0}$ and $\{f_j^{-1}\circ f_j=\id_{A^j}\}_{j\geq 0}$ both have inverse limits.  Let $g\colon A^\omega\to A^\omega$ be the inverse limit of the former.  The latter inverse limit is $\id\colon A^\omega\to A^\omega$, but by uniqueness, $\id=g\circ f$.  Similarly, $\id=f\circ g$, which shows $f$ is bijective.  
\end{proof}

\begin{example}[Functions on $\Z_2$]
In order to discuss concrete examples, we will often fix $A=2$, the
two-element set $\{0, 1\}$.  The set of bitstreams $2^\omega$
underlies the ring $\Z_2$ of $2$-adic integers.  The following
examples show how $2$-adic arithmetic operations correspond to the
notion of $k$-causal functions on the underlying streams.
\begin{enumerate}[(i)]
\item The mappings $x\mapsto 2x$ and $x\mapsto 2x+1$, respectively,
correspond to $c_0$ and $c_1$ on $2^\omega$, and consequently, they
are both $1$-causal.
\item The function $t\colon\Z_2\to\Z_2$ given by
\begin{equation}
t(x) = 
\begin{cases}
\frac{x}{2} & \text{if $x\equiv 0\pmod{2}$}\\
\frac{x-1}{2} & \text{if $x\equiv 1\pmod{2}$}
\end{cases}
\label{eqn:tail-def2}
\end{equation}
corresponds exactly to the tail map on streams (\ref{eqn:tail-def})
and is $(-1)$-causal.  In the case of $2$-adic integers
(i.e.~bitstreams), we call $x$ even or odd, depending on whether
$h(x)=0$ or $h(x)=1$. The two cases in (\ref{eqn:tail-def2})
differentiate between the two possibilities.

\item For any $k\in\Z_2$, the mapping $x\mapsto x+ k$ is bicausal.
For any $k\in\Z_2$, the mapping $x\mapsto k\cdot x$ is causal, but
only bicausal if $k$ is a unit (via Corollary
\ref{cor:bicausal-finite-alphabet}(\ref{item:causal-bijection})).  Since the composition
of bicausals is bicausal, the mapping $x\mapsto ax+b$ is bicausal for
any $b\in\Z_2$ and $a\in 2\Z_2+1$.
\end{enumerate}
\end{example}

\section{Coalgebras and Coinduction}
\label{sec:coalgebra}

In a very general sense, coinduction is a notion which is
categorically dual to induction.  Here we will introduce the bare
minimum of ideas in the theory of coalgebras, but more thorough
introductions are available elsewhere \cite{Ru2}.  Though we start
with the most general (category theoretic) definition, we will quickly focus
on a particular instance of stream coalgebras.

\begin{defn}
Given an endofunctor $F$ on a category $\Ccal$, an $F$-coalgebra is a
$\Ccal$-morphism $X\xrightarrow{f}FX$.  An \emph{$F$-coalgebra
morphism} from $X\xrightarrow{f}FX$ to $Y\xrightarrow{g}FY$ is a
$\Ccal$-morphism $X\xrightarrow{m}Y$ such that the diagram
$$
\xymatrix{
X \ar[rr]^{f} \ar[d]_{m} & & FX \ar[d]^{Fm}\\
Y \ar[rr]^{g} & & FY
}
$$ commutes (i.e.~$g\circ m = Fm\circ f$). The class of all
$F$-coalgebras and $F$-coalgebra morphisms forms a category
$\Coalg(F)$. Terminal objects in $\Coalg(F)$, if they exist at all,
are called \emph{final $F$-coalgebras}.  If $\Ccal$ is the category
$\Set$ of sets---as assumed in the rest of this paper---we often refer
to the coalgebra as the \emph{structure map} and its domain as the
\emph{carrier set}.
\end{defn}

\begin{example}[Trivial]
One of the most trivial examples is the coalgebra where the
endofunctor is the identity functor.  Coalgebras in this context are
endofunctions where the domain and codomain coincide.  A coalgebra
morphism from $\left(X\xto{f} X\right)$ to $\left(Y\xto{g} Y\right)$
is a function $X\xto{m}Y$ so that $g\circ m = m\circ f$.  Any
singleton set with the identity function is a final coalgebra.
\end{example}

\begin{example}[Mealy automaton]
For a more substantial example, consider the endofunctor $M_{A,B}$
given by $M_{A,B}X=(B\times X)^A$ for a pair of sets $A$ and $B$.
While this example will only be incidental to the main results of this
paper, it is interesting to note that $0$-causal functions appear
naturally in other coalgebraic situations.  Coalgebras in this context
have the form $X\xrightarrow{\alpha}(B\times X)^A$ and are called
\emph{Mealy automata} with input and output taken from sets $A$ and
$B$, respectively. The coalgebra structure map $\alpha$ corresponds to
a (deterministic) transition function $A\times X\xrightarrow{\check
\alpha}B\times X$ on the state space $X$.  Here $\check\alpha$ gives
an output symbol from $B$ and the next state in $X$ given an input
symbol from $A$ and the current state in $X$.

Rutten \cite{Ru1} showed that the set $\Gamma_0$ of all $0$-causal
functions from $A^\omega$ to $B^\omega$ forms the carrier set of a
final coalgebra $\Gamma_0\xrightarrow{\gamma} (B\times \Gamma_0)^A$
for the Mealy automaton endofunctor $M_{A,B}$.  The structure map $\gamma$ is given
by
$$
\gamma(f)(a) =  (D_1(f\circ c_a), t\circ f\circ c_a)
$$ for $f\in\Gamma_0$, $a\in A$, and $\sigma\in A^\omega$.
The finality in $\Coalg(M_{A,B})$ amounts to the idea
that given any Mealy automaton, each state $x_0$ can be assigned a
$0$-causal function which encodes the output of the automaton for any
stream of inputs that might be presented to the automaton starting at
state $x_0$.  This illustrates a theme in the theory of coalgebras
that morphisms into final coalgebras often encapsulate infinite
dynamics of any given coalgebra.
\end{example}

\begin{example}[Power set functor]
For a more pathological example, consider the functor $P$ given
$PX=\mathcal{P}(X)$ (the power set of $X$).  $P$-coalgebras are
exactly graphs, encoded as functions of the type $X\to PX$.  Cantor's
theorem tells us that $PX$ has a strictly greater cardinality than
$X$.  However, Lambek's Lemma asserts that final coalgebras must be
isomorphisms (i.e.~bijections in $\Set$).  This proves that
$\Coalg(P)$ has no final coalgebra.
\end{example}

For the rest of the paper, we will focus our attention exclusively on
stream coalgebras (over $\Set$).  In the case where the endofunctor
$S_B$ is given by $S_B X=B\times X$ for set $B$, we call these
coalgebras \emph{$B$-stream coalgebras}.  Regardless of what $B$ might
be, final stream coalgebras do exist, and the standard example is
$$
B^\omega\xrightarrow{\tuple{h,t}} B\times B^\omega
$$ 
where $h\colon B^\omega\to B$ and $t\colon B^\omega\to B^\omega$ are the \emph{head} and \emph{tail} maps defined earlier.  Finality for $\Coalg(S_B)$ can be stated in the
following way.  Given any $B$-stream coalgebra,
$X\xrightarrow{\tuple{g,s}} B\times X$, there exists a unique map
$\Phi=\Phi_{\tuple{g,s}}\colon X\to B^\omega$ so that the diagram
\begin{equation}
\label{eqn:induced-homomorphism}
\xymatrix{
X \ar[rr]^{\tuple{g,s}} \ar[d]_{\Phi} & & B\times X \ar[d]^{\id_B\times \Phi} \\
B^\omega \ar[rr]_{\tuple{h,t}} & & B\times B^\omega
}
\end{equation}
commutes, i.e.~$g=h\circ \Phi$ and $t\circ \Phi = \Phi\circ s$.  We
say the unique map $\Phi$ above is \emph{coinductively induced} (or
equivalently, \emph{coiteratively defined}) by $\tuple{g,s}$.  We can
actually define $\Phi$ from $\tuple{g,s}$ explicitly:
\begin{equation}
\Phi(x)(n) = h(t^{(n)}(x))= g(s^{(n)}(x))
\label{eqn:stream-coinduced}
\end{equation}
As a stream, $\Phi(x)$ records the image of the $s$-orbit of $x$
under $g$.

\section{Results}
\label{sec:results}

For the main theorem, we consider the properties of a function $\phi$ coinductively induced by a coalgebra $A^\omega\xto{\tuple{H,T}} B\times A^\omega$ whose carrier set is the set of streams and whose structure map incorporates a function $T$ woven from a set of $0$-causal functions.  The resulting map $\phi$ is therefore a stream function, and we explore its relationship to the set of stream function from which it was derived.

\begin{thm}
Let $A^\omega\xto{\tuple{H,T}} B\times A^\omega$ be a $B$-stream coalgebra so that 
\begin{equation}
h(\sigma)=h(\tau) \quad\Longrightarrow\quad
H(\sigma)=H(\tau)
\label{eqn:thm}
\end{equation}
for all $\sigma,\tau\in A^\omega$.
Let $\phi\colon
A^\omega\to B^\omega$ be the (unique) coalgebra morphism induced
by $\tuple{H,T}$.
\begin{enumerate}[(i)]
\item The coalgebra morphism $\phi$ is
causal if $T$ is woven from a family of $0$-causal functions. 
\item If the converse of (\ref{eqn:thm}) also holds, then the
coalgebra morphism $\phi$ is bicausal if and only if $T$ is woven from
a family of bicausal functions.
\end{enumerate}
\label{thm:main}
\end{thm}

\begin{remark}
Because $\phi$ is coinductively induced by $\tuple{H,T}$, we have two identities: $h\circ\phi = H$ and $t\circ\phi = \phi\circ T$ (or more generally, $t^{(n)}\circ\phi = \phi\circ T^{(n)}$ for all $n\geq 0$).  The proof uses the second characterization of $0$-causal and bicausal functions from Corollaries \ref{cor:0-causal} and \ref{cor:bicausal}, respectively.  For each $n\geq 0$, $0$-causal and bicausal functions must preserve $n$-prefix-equivalence.  The proofs proceed by induction on $n\geq 0$.  
\end{remark}

\begin{proof}
The base case is trivial since $\sigma\equiv_0 \tau$ for all $\sigma,\tau\in A^\omega$.  For the induction case, assume for some $n\geq 0$,
\begin{equation}
\sigma\equiv_n\tau \quad\Longrightarrow\quad \phi(\sigma)\equiv_n \phi(\tau)
\end{equation}
and suppose $\sigma\equiv_{n+1}\tau$.  First of all, since $T$ is woven from $0$-causal functions, it is $(-1)$-causal, and therefore $T^{(n)}$ is $(-n)$-causal.  By Theorem \ref{thm:k-causal-equivalence}, ${T^{(n)}(\sigma)\equiv_1 T^{(n)}(\tau)}$, or equivalently, $h(T^{(n)}(\sigma))=h(T^{(n)}(\tau))$.  The premise (\ref{eqn:thm}) of the theorem therefore gives us: 
\begin{equation}
H(T^{(n)}(\sigma)) = H(T^{(n)}(\tau)).
\label{eqn:middle-step}
\end{equation}  
Secondly, $\sigma\equiv_{n+1}\tau$ implies that $\sigma\equiv_n\tau$.  By the induction hypothesis, $\phi(\sigma)\equiv_n \phi(\tau)$.   Therefore, to check that $\phi(\sigma)\equiv_{n+1} \phi(\tau)$, it is only necessary to verify that $\phi(\sigma)$ and $\phi(\tau)$ agree at index $n$:
\begin{align*}
\phi(\sigma)(n) 
&= h(t^{(n)}(\phi(\sigma))) &&[\text{Lemma \ref{lem:observation}(\ref{lem:observation-4})}] \\
&= h(\phi(T^{(n)}(\sigma))) &&[t\circ\phi = \phi\circ T]\\
&= H(T^{(n)}(\sigma))  &&[h\circ\phi = H]\\
&= H(T^{(n)}(\tau))  &&[\text{(\ref{eqn:middle-step})}]\\
&= h(\phi(T^{(n)}(\tau)))
= h(t^{(n)}(\phi(\tau)))
= \phi(\tau)(n)
\end{align*}
Therefore $\phi(\sigma)\equiv_{n+1}\phi(\tau)$.  This completes the proof of (i).

For the bicausal case, suppose $T$ is woven from bicausal functions.  The basic step of the induction is the same as above.  For the induction step, assume for some $n\geq 0$,
\begin{equation}
\sigma\equiv_n\tau \Longleftrightarrow \phi(\sigma)\equiv_n \phi(\tau)
\label{eqn:ind-hyp-2}
\end{equation}
We have already showed above that this implies $\sigma\equiv_{n+1}\tau \Longrightarrow \phi(\sigma)\equiv_{n+1} \phi(\tau)$, so consider the converse, and suppose $\phi(\sigma)\equiv_{n+1}\phi(\tau)$.  In particular, $\phi(\sigma)\equiv_1\phi(\tau)$, i.e.~$h(\phi(\sigma)) = h(\phi(\tau))$.  Since $h\circ\phi=H$, we have $H(\sigma)=H(\tau)$.  By the premise of the theorem---the converse of (\ref{eqn:thm}), to be precise---we can conclude that $h(\sigma)=h(\tau)$.  Let $a=h(\sigma)=h(\tau)$, and proceed:
\begin{align*}
\phi(\sigma)\equiv_{n+1}\phi(\tau)
\quad\Rightarrow\quad & t(\phi(\sigma))\equiv_n t(\phi(\tau)) & &\text{[$t$ is $(-1)$-causal]}\\
\quad\Leftrightarrow\quad & \phi(T(\sigma))\equiv_n \phi(T(\tau)) & &\text{[$t\circ\phi = \phi\circ T$]}\\
\quad\Leftrightarrow\quad & T(\sigma)\equiv_n T(\tau) & &\text{[induction hypothesis (\ref{eqn:ind-hyp-2})]}\\
\quad\Leftrightarrow\quad & f_{h(\sigma)}(t(\sigma))\equiv_n f_{h(\tau)}(t(\tau)) & &\text{[$T$ is woven from $\{f_\alpha\}$]}\\
\quad\Leftrightarrow\quad & t(\sigma)\equiv_n t(\tau) & &\text{[$f_\alpha$ is bicausal]}\\
\quad\Leftrightarrow\quad & a\append t(\sigma)\equiv_{n+1} a\append t(\tau) & &\text{[$c_a$ is $1$-causal]}\\
\quad\Leftrightarrow\quad & \sigma\equiv_{n+1}\tau& &\text{[$a=h(\sigma)=h(\tau)$]}
\end{align*}

For the other direction of (ii), suppose $\phi$ is bicausal.  Then, for $\sigma,\tau\in A^\omega$,
\begin{align*}
d(T(\sigma),T(\tau)) 
&= d(\phi(T(\sigma)), \phi(T(\tau))) &&[\text{$\phi$ is bicausal}]\\
&= d(t(\phi(\sigma)),t(\phi(\tau))) &&[t\circ \phi = \phi\circ T]\\
&\leq 2d(\phi(\sigma),\phi(\tau)) &&[\text{$t$ is $(-1)$-causal}]\\
&= 2d(\sigma,\tau) &&[\text{$\phi$ is bicausal}]
\end{align*}
This shows that $T$ is $(-1)$-causal and therefore woven from $0$-causal functions, by Lemma
\ref{lem:k-causal-is-woven}.  Recall from the proof of the lemma that
$T$ is woven from $\{f_a\}_{a\in A}$ where each $f_a$ is given by
$f_a(\sigma) = T(a\append\sigma)$.  We must show that $f_a$ is bicausal for any $a\in A$, but first recall that since $\phi$ is $0$-causal, the map $\phi\circ c_a$ is $1$-causal.  In particular, the head of $\phi(c_a(\sigma))$ does not depend on $\sigma$.  Therefore we have 
\begin{equation}
t(\phi(c_a(\sigma)))\equiv_{j}t(\phi(c_a(\tau)))
\quad\Longleftrightarrow\quad
\phi(c_a(\sigma))\equiv_{j+1}\phi(c_a(\tau)) 
\label{eqn:awkward}
\end{equation}  
for all $\sigma,\tau\in A^\omega$, $j\in\N$. Then we can proceed:
\begin{align*}
f_a(\sigma)\equiv_n f_a(\tau)
\quad\Leftrightarrow\quad & T(a\append\sigma)\equiv_n T(a\append\tau) & &\text{[def. of $f_a$]}\\
\quad\Leftrightarrow\quad & \phi(T(a\append\sigma))\equiv_n \phi(T(a\append\tau)) & &\text{[$\phi$ bicausal]}\\
\quad\Leftrightarrow\quad & t(\phi(a\append\sigma))\equiv_n t(\phi(a\append\tau)) & &\text{[$t\circ\phi = \phi\circ T$]}\\
\quad\Leftrightarrow\quad & \phi(a\append\sigma)\equiv_{n+1} \phi(a\append\tau) & &\text{[(\ref{eqn:awkward})]}\\
\quad\Leftrightarrow\quad & a\append\sigma\equiv_{n+1} a\append\tau & &\text{[$\phi$ bicausal]}\\
\quad\Leftrightarrow\quad & \sigma \equiv_n \tau & &\text{[$t$ is $(-1)$-causal]}
\end{align*}
This calculation shows that $f_a$ is distance-preserving.
\end{proof}



By taking $A=B$ and $H=h$, we can immediately get the following corollary.

\begin{corollary}
Let $T\colon A^\omega\to A^\omega$ be a function.  Let $\phi\colon
A^\omega\to A^\omega$ be the (unique) coalgebra morphism induced by
the coalgebra $A^\omega\xrightarrow{\tuple{h,T}} A\times A^\omega$.
\begin{enumerate}[(i)]
\item The coalgebra morphism $\phi$ is
causal if $T$ is woven from a family of $0$-causal functions. 
\item The coalgebra morphism $\phi$ is bicausal if and only if $T$ is
woven from a family of bicausal functions.
\end{enumerate}
\label{cor:woven-causal-induce-causals}
\end{corollary}

In essence, this corollary asserts that the sets $\Gamma_0$ of
$0$-causal functions and $\Gammabi$ of bicausal functions on a set of
streams are both closed under a particular coinductive construction.
If $T\colon A^\omega\to A^\omega$ is woven from a subset of $\Gamma_0$
(resp. $\Gammabi$), then the coalgebra morphism $\phi$ coinductively
induced by $\tuple{h,T}$ lies in $\Gamma_0$ (resp. $\Gammabi)$.

Theorem \ref{thm:main} and Corollary
\ref{cor:woven-causal-induce-causals} generalize a result for $A=B=2$
and ``solenoidal bijections'' in a paper by Bernstein and Lagarias
\cite{BL} in several ways.  First, it is useful to know that the ring
structure of $\Z_2$ is not strictly necessary.  Second, the
cardinality of the underlying set is not required to be two, prime, or
even finite.  Lastly, if one direction of the implication is
abandoned, instead of requiring bicausal functions, the premise of the theorem can be weakened to
$0$-causal functions.

\subsection*{Connection to the $3x+1$ Problem}

In the $2$-adic context, if we specify that $f_0(x)=x$ and
$f_1(x)=3x+2$, the resulting woven function is
\begin{equation}
T(x) = 
\begin{cases}
\frac{x}{2} & \text{if $h(x)=0$}\\
3\frac{x-1}{2}+2 & \text{if $h(x)=1$}
\end{cases}
\label{eqn:collatz}
\end{equation}
Noting that $3\frac{x-1}{2}+2 = \frac{3x+1}{2}$, we can see that
(\ref{eqn:collatz}) is the definition of the $2$-adic extension of
Collatz function, which on the integers is given by: 
$$
C(n) = 
\begin{cases}
\frac{n}{2} & \text{if $n$ even}\\
\frac{3n+1}{2} & \text{if $n$ odd}
\end{cases}
$$ The famed $3x+1$ Problem is to determine whether or not for all
$n>0$, there exists a $k$ where $C^{(k)}(n)=1$.  All computational
evidence point to an affirmative answer; as of June 2008 \cite{OS1},
the conjecture has been verified up to $n=18\cdot 2^{58}$ by machine,
but the problem is still unsolved in general \cite{Lag1,Lag2,Ak}.

Noting that $f_0$ and $f_1$ are both bicausal, we conclude that the
coalgebra morphism $Q\colon \Z_2\to 2^\omega$ coinductively induced by
$\tuple{h,T}$ is bicausal.  Since $A=2$ is finite, $Q$ is a bijection.
Moreover, bijective coalgebra morphisms are also coalgebra
isomorphisms \cite{Ru1}. In other words, we have the following result.

\begin{thm}
Let $T$ be the $2$-adic extension of the $3x+1$ function given in
(\ref{eqn:collatz}).  In the category of 2-stream coalgebras,
$\tuple{h,T}$ is terminal.  More specifically, for any 2-stream
coalgebra $X\xto{f}2\times X$, there exists a unique coalgebra
morphism $X\xto{\psi} \Z_2$ so that $h\circ \psi = h$ and $T\circ \psi = \psi\circ f$.
\label{thm:collatz-final}
\end{thm}

\begin{remark}
The unique coalgebra isomorphism $Q$ here is the \emph{$3x+1$
conjugacy map} denoted in Lagarias \cite{Lag0} as $Q_\infty$.  It
produces the \emph{parity vector} for $2$-adic integers,
i.e.~$Q(x)(n)$ is the parity of the $n$th iterate of $x\in\Z_2$ under
$T$.  It should be clear that this is a special case of
(\ref{eqn:stream-coinduced}).  In addition, Bernstein \cite{Be}
exhibited a formula for the inverse of $Q$ (using the notation from
this paper):
$$
Q^{-1}(\sigma) = -\frac{1}{3} \sum_{\sigma(\ell)=1} \frac{1}{3^\ell} 2^\ell .
$$

\end{remark}

This coinductive observation about $T$ also
applies to any function woven from bicausals on streams with a finite
alphabet.  For instance, for $m,n\in \Z$, let $T_{m,n}\colon \Z_2\to
\Z_2$ be given by
\begin{equation}
T_{m,n}(x) = 
\begin{cases}
Q^{(m)}(\frac{x}{2})& \text{if $h(x)=0$}\\
Q^{(n)}(\frac{x-1}{2}) & \text{if $h(x)=1$}
\end{cases}
\end{equation}
woven from $f_0=Q^{(m)}$ and $f_1=Q^{(n)}$.  By the same argument as
for $T$, we can conclude that $\tuple{h,T_{m,n}}$ is a final stream
coalgebra.  Given the plethora of other examples, coalgebraic finality
in and of itself cannot yield a solution to the $3x+1$
problem. Nonetheless, it is worth noting that many of the interesting
cousins of the $3x+1$ Problem are based on the dynamics of woven
functions (often woven from bicausals).

\end{document}